\documentclass{article}
\usepackage{latexsym,amssymb}
\input amssym.def \input amssym

\newtheorem{te}{Theorem}[section]

\newtheorem{de}[te]{Definition}
\newtheorem{lm}[te]{Lemma}
\newtheorem{fa}[te]{Fact}

\newtheorem{ex}[te]{Example}

\newcommand{\fun}[2]{{{}^{#1}\hspace{-1mm}#2}}
\def\dokaz{\noindent{\bf Proof. }}
\def\kraj{\hfill $\Box$ \par \vspace*{2mm} }
\def\widemid{\hspace{1mm}\widetilde{\mid}\hspace{1mm}}
\def\nwidemid{\hspace{1mm}\widetilde{\nmid}\hspace{1mm}}

\def\po{\exists}
\def\str{\rightarrow}

\def\dl{\Leftrightarrow}
\def\dstr{\hspace{-0.1cm}\downarrow}
\def\gstr{\hspace{-0.1cm}\uparrow}
\def\ps{\subseteq}

\def\ra{\setminus}

\def\rest{\upharpoonright}

\def\cU{{\cal U}}
\def\cV{{\cal V}}
\def\cA{{\cal A}}
\def\cB{{\cal B}}

\def\cF{{\cal F}}

\begin{document}
\begin{center}
           {\large \bf $\widemid$-divisibility of ultrafilters}
\end{center}
\begin{center}
{\small \bf Boris \v Sobot}\\[2mm]
{\small Faculty of Sciences, University of Novi Sad,\\
Trg Dositeja Obradovi\'ca 4, 21000 Novi Sad, Serbia\\
e-mail: sobot@dmi.uns.ac.rs}
\end{center}
\begin{abstract} \noindent
We further investigate a divisibility relation on the set $\beta N$ of ultrafilters on the set of natural numbers. We single out prime ultrafilters (divisible only by 1 and themselves) and establish a hierarchy in which a position of every ultrafilter depends on the set of prime ultrafilters it is divisible by. We also construct ultrafilters with many immediate successors in this hierarchy and find positions of products of ultrafilters.
\vspace{1mm}\\


\end{abstract}

\section{Introduction}

Let $N$ denote the set of natural numbers (without zero). The Stone-\v Cech compactification of the discrete space $N$ is the space $\beta N$ of all ultrafilters over $N$. For each $n\in N$ the principal ultrafilter $p_n=\{A\ps N:n\in A\}$ is identified with the respective element $n$. The topology on $\beta N$ is generated by (clopen) base sets of the form ${\bar A}=\{x\in\beta N:A\in x\}$.

A family $\cF$ of subsets of $N$ has the finite intersection property (f.i.p.) if the intersection of every finitely many elements of $\cF$ is nonempty. $\cF$ has the uniform f.i.p.\ if the intersection of every finitely many elements of $\cF$ is infinite. Every family with f.i.p.\ is contained in an ultrafilter, and every family with u.f.i.p.\ is contained in a nonprincipal ultrafilter. A family $\cF$ with the f.i.p.\ generates a filter $F$ if for every $B\in F$ there are $A_1,A_2,\dots,A_n\in\cF$ such that $A_1\cap A_2\cap\dots\cap A_n\subseteq B$.

If $f:N\str N$ is a function, the direct and inverse image of a set $A\subseteq N$ are $f[A]=\{f(a):a\in A\}$ and $f^{-1}[A]=\{b:f(b)\in A\}$. For every $f:N\str N$ there is unique continuous function $\widetilde{f}:\beta N\str\beta N$ extending $f$. It is given by $\widetilde{f}(x)=\{A\ps N:f^{-1}[A]\in x\}$ for $x\in\beta N$. $\widetilde{f}(x)$ is also generated by sets $f[A]$ for $A\in x$. Since $\widetilde{g}\circ\widetilde{f}$ is a continuous extension of $g\circ f$, it follows that $\widetilde{g}\circ\widetilde{f}=\widetilde{g\circ f}$ for every two functions $f,g:N\str N$.

The multiplication $\cdot$ on $N$ can be extended to $\beta N$ (using the same notation $\cdot$ for the extension) in such way that $(\beta N,\cdot)$ is a compact Hausdorff right-topological semigroup: for $x,y\in\beta N$,
\begin{equation}\label{product}
A\in x\cdot y\dl\{n\in N:A/n\in y\}\in x,
\end{equation}
where for $A\ps N$ and $n\in N$, $A/n=\{m\in N:mn\in A\}=\{\frac an:a\in A,n\mid a\}$. The known properties if this structure are described in detail in \cite{HS}.

We fix some more notation. Throughout the paper $P$ is the set of prime numbers. The complement of $A\subseteq N$ with respect to $N$ is $A^c=N\setminus A$, and the complement of $A\subseteq P$ with respect to $P$ is $A'=P\setminus A$. $[X]^k$ is the set of subsets of $X$ of cardinality $k$. The set of functions $f:A\str B$ will be denoted by $\fun AB$.

For $x\in\beta N$, if $A\in x$ then $X\in x$ if and only if $X\cap A\in x$. Because of this we can identify each ultrafilter $x$ on $N$ containing $A$ with the ultrafilter $x\rest A=\{X\cap A:X\in x\}$ on $A$. Thus it is also common to think of $\overline{A}$ as a subspace of $\beta N$. Also, $A^*=\overline{A}\ra A$ is the set of nonprincipal ultrafilters containing $A$.

The Rudin-Keisler preorder on $\beta N$ is defined as follows: $x\leq_{RK} y$ if and only of there is $f:N\str N$ such that $\widetilde{f}(y)=x$.

An ultrafilter $x\in N^*$ is a P-point if, for every sequence $\langle A_n:n\in N\rangle$ of sets in $x$ such that $A_m\supseteq A_n$ for $m<n$ there is a set $B\in x$ (called a pseudointersection of sets $A_n$) such that $B\setminus A_n$ is finite for every $n\in N$.

A function $c:X\str\{0,1\}$ is called a $2$-coloring of $X$. An ultrafilter $x$ is called Ramsey (selective) if, for every $A\in x$, every $k\in N$ and every $2$-coloring $c$ of $[A]^k$ there is a set $M\in x$ that is monochromatic, i.e.\ such that $|c[[M]^k]|=1$. It is well-known that $x$ is a Ramsey ultrafilter if and only if $x$ is minimal in $\leq_{RK}$ (there are no nonprincipal ultrafilters $y\leq_{RK}x$).

If $A\subseteq N$, we will denote $A\gstr=\{n\in N:\po a\in A\;a\mid n\}$ and $A\dstr=\{n\in N:\po a\in A\;n\mid a\}$. $\cU=\{A\ps N:A=A\gstr\}$ and $\cV=\{A\ps N:A=A\dstr\}$ are the collections of subsets of $N$ upwards/downwards closed for divisibility.

In \cite{So1} four divisibility relations on $\beta N$ were introduced. In \cite{So1} and \cite{So2} some properties of these relations were investigated. In particular, the relation $\widemid$ was introduced as an extension of the usual divisibility relation $\mid$ on $N$ to $\beta N$ analogous to extensions of functions described above. It was proven that, for $x,y\in\beta N$,
\begin{eqnarray*}
x\widemid y &\dl& x\cap\cU\subseteq y\\
 &\dl& y\cap\cV\subseteq x,
\end{eqnarray*}
In this paper we will use these characterizations and under divisibility of ultrafilters we will understand $\widemid$-divisibility. For example, "least upper bound for $x$ and $y$" will mean the $\widemid$-smallest ultrafilter $z$ (if such exists) such that $x\widemid z$ and $y\widemid z$. "$x$ is below $y$" will mean $x\widemid y$. $\widemid$ is a preorder and, if we define $x=_\sim y\Leftrightarrow x\widemid y\land y\widemid x$, then $=_\sim$ is an equivalence relation. We denote the equivalence class of $x$ by $[x]_\sim$ and think of $\widemid$ as an order on the set of these classes. Let us also call $x,y\in\beta N$ incompatible if there is no $z\in\beta N\setminus\{1\}$ such that $z\widemid x$ and $z\widemid y$.

The motivation for introducing divisibility of ultrafilters is to inspect the effects of existence of some infinite sets of natural numbers to ultrafilters and (hopefully) the opposite: to draw number-theoretic conclusions from the existence of certain ultrafilters. However, we also hope that some of the results in this paper will help to better understand the product (\ref{product}) of ultrafilters by finding its place in the divisibility hierarchy.

\section{Prime ultrafilters}\label{secwidemin}

\begin{lm}\label{funkcija}
Let $x\in\beta N$, $A\in x$ and $f:N\str N$.

(a) If $f(a)\mid a$ for all $a\in A$, then $\widetilde{f}(x)\widemid x$.

(b) If $a\mid f(a)$ for all $a\in A$, then $x\widemid \widetilde{f}(x)$.
\end{lm}

\dokaz (a) If $B\in{\widetilde f}(x)\cap\cU$ then $f^{-1}[B]\in x$. But $f^{-1}[B]\cap A\subseteq B$ (because $B\in\cU$) so $B\in x$. Hence we have ${\widetilde f}(x)\widemid x$.

(b) Analogously to (a) we prove that every $B\in{\widetilde f}(x)\cap\cV$ is also in $x$.\kraj

If $A\in x$, in order to determine $\widetilde{f}(x)$ it is enough to know values of $f(a)$ for $a\in A$. Hence, when using the lemma above we will sometimes define functions only on a set in $x$.

\begin{lm}\label{minbelow}
For every $x\in\beta N\setminus\{1\}$ there is $p\in\overline{P}$ such that $p\widemid x$.
\end{lm}

\dokaz We define a function $f:N\setminus\{1\}\str N$: let $f(n)$ be the smallest prime factor of $n$. Now ${\widetilde f}(x)\in\overline{P}$ (for $x\in\beta N\setminus\{1\}$) and, by Lemma \ref{funkcija}(a), ${\widetilde f}(x)\widemid x$.\kraj

Clearly, 1 is the smallest element in $(\beta N,\widemid)$. We will call $p\in\beta N\setminus\{1\}$ {\it prime} (or $\widemid$-minimal) if it is divisible only by 1 and itself. We will reserve labels $p,q,r,\dots$ for prime ultrafilters and $x,y,z,\dots$ for ultrafilters in general.

\begin{te}\label{widemin}
$p\in\beta N$ is prime if and only if $p\in\overline{P}$.
\end{te}

\dokaz Assume first that $p$ is prime but $p\in\overline{P^c}$. By Lemma \ref{minbelow} there is an element $q\in\overline{P}$ below $p$. But $q\neq_\sim p$ because the set of composite numbers $(P\cup\{1\})^c\in(p\cap\cU)\setminus q$; hence $q\neq p$, a contradiction.

Now assume $p\in\overline{P}$ and $x\widemid p$, but $x\neq p$. By Lemma \ref{minbelow} there is $q\in\overline{P}$ such that $q\widemid x$. It follows that $q\widemid p$. If $A\in q\setminus p$, let $B=(A\cap P)\uparrow$. Then $B\in q\cap\cU$, so $B\in p$ and $A\supseteq B\cap P\in p$, a contradiction.\kraj

Comparing this with Lemmas 7.3 and 7.5 from \cite{So1} we conclude that $\widemid$-minimal ultrafilters are also $\mid_L$-minimal, but not vice versa. Also, a corollary of this theorem is that there is a family of $2^{\goth c}$ incompatible ultrafilters, which improves Theorem 3.8 from \cite{So2}.

Let us also mention that, by Fact \ref{greatest}, there are ultrafilters with $2^{\goth c}$-many divisors so not every divisibility can be established by means of Lemma \ref{funkcija}(a). It follows that the Rudin-Keisler preorder is not stronger than $\widemid$. On the other hand, neither is $\widemid$ stronger than $\leq_{RK}$, since ultrafilters containing $P$ that are not Ramsey are $\widemid$-minimal but nor $\leq_{RK}$-minimal.

Now we define levels of the $\widemid$-hierarchy.

\begin{de}
Let $L_0=\{1\}$ and $L_n=\{a_1a_2\dots a_n:a_1,a_2,\dots,a_n\in P\}$ for $n\geq 1$. We will say that $x$ is of level $n$ if $x\in\overline{L_n}$.
\end{de}

In particular, $\overline{L_1}=\overline{P}$ is the set of prime ultrafilters.

In this paper we will mostly deal with ultrafilters on finite levels (belonging to $\overline{L_n}$ for some $n\in N$). In Section \ref{above} we will see that there are also ultrafilters above all these finite levels.

\section{Ultrafilters with one prime divisor}

\begin{de}\label{defstepen}
For $A\subseteq N$ and $n\in N\cup\{0\}$ we denote $A^n=\{a^n:a\in A\}$. If $pow_n:N\str N$ is defined by $pow_n(a)=a^n$ then, for $x\in N^*$, $\widetilde{pow_n}(x)$ is generated by sets $A^n$ for $A\in x$. We will denote $\widetilde{pow_n}(x)$ with $x^n$.
\end{de}

Of course, $A^0=\{1\}$, $x^0=1$ and $x^1=x$.

\begin{lm}\label{kvadrati}
If $p\in\overline{P}$, the only ultrafilters below $p^n$ are $p^k$ for $k\leq n$.
\end{lm}

\dokaz If $p\in P$, the lemma is obvious. So we prove it for $p\in P^*$.

Since $\bigcup_{k\leq n}P^k\in p^n\cap\cV$, every ultrafilter $x\neq 1$ below $p^n$ must contain $P^k$ for some $k\leq n$.

If $x\in(P^k)^*$ then for every $A\subseteq P$, if $A^k\in p^k$ then we have $\bigcup_{i\leq n}A^i\in p^n\cap\cV$ so $x\widemid p^n$ implies that $A^k=P^k\cap\bigcup_{i\leq n}A^i\in x$. Thus $x\rest P^k=p^k\rest P^k$, which means that $x=p^k$.\kraj

\begin{de}
For $A,B\subseteq N$ let $AB=\{ab:a\in A,b\in B\land gcd(a,b)=1\}$. In particular, $A^{(n)}=\underbrace{A\cdot A\cdot\dots\cdot A}_n$. If $p\in\overline{P}$ let $F_n^p=\{A^{(n)}:A\in p\rest P\}$.
\end{de}

Note the difference between $A^{(n)}$ and $A^n$ from Definition \ref{defstepen}: elements of $A^{(n)}$ must be products of mutually prime numbers. $A$ will almost always be a subset of $P$ in which case "mutually prime" will mean "distinct".

If $p\in P$ (a principal ultrafilter), $\{p\}^{(n)}=\emptyset\in F_n^p$, so there is no ultrafilter $x\supseteq F_n^p$.

\begin{lm}
Let $p\in P^*$. For any ultrafilter $x\supseteq F_n^p$, $p$ is the only ultrafilter from $\overline{P}$ below $x$.
\end{lm}

\dokaz If $q\in\overline{P}$ is such that $q\neq p$, there is $A\subseteq P$ such that $A\in p\setminus q$. Then $\bigcup_{k\leq n}A^{(k)}\in x\cap\cV\setminus q$, so $q\widetilde{\nmid}x$.\kraj

\begin{lm}\label{broj}
Let $p\in P^*$.

(a) $p\cdot p\supseteq F_2^p$.

(b) There are either finitely many or $2^{\goth c}$ ultrafilters $x\supseteq F_n^p$.
\end{lm}

\dokaz (a) Let $A\in p\rest P$. Since $A^{(2)}/a=A\setminus\{a\}$ for $a\in A$, $\{a\in N:A^{(2)}/a\in p\}\supseteq A\in p$ and so $A^{(2)}\in p\cdot p$.\\

(b) The set of ultrafilters $x\supseteq F_n^p$ is actually $\bigcap_{X\in F_n^p}\overline{X}$ so it is closed. By \cite{W}, Theorem 3.3, closed subsets of $\beta N$ are either finite or of cardinality $2^{\goth c}$. Hence, for every $p\in P^*$ there are either finitely many or $2^{\goth c}$ ultrafilters containing $F_n^p$.\kraj

\begin{te}\label{ramsey}
Let $p\in P^*$. There is unique ultrafilter $x\supseteq F_n^p$ if and only if $p$ is Ramsey.
\end{te}

\dokaz Assume $p$ is Ramsey. To prove that $x$ is unique it suffices to show that for every set $S\subseteq P^{(n)}$ one of the sets $S$ and $P^{(n)}\setminus S$ contains a set from $F_n^p$. Define a coloring of $[P]^n$ as follows:
$$c(\{a_1,a_2,\dots,a_n\})=\left\{\begin{array}{ll}
0, \mbox{if }a_1a_2\dots a_n\in S\\
1, \mbox{otherwise.}
\end{array}\right.$$
Since $p$ is Ramsey there is a monochromatic set $M\in p$, say $c(\{a_1,a_2,\dots,a_n\})=0$ for all $a_1,a_2,\dots,a_n\in M$. This means that $M^{(n)}\subseteq S$.

Now assume $p$ is not Ramsey. Then there is a 2-coloring $c$ of $[P]^n$ such that $p$ does not contain a monochromatic subset. Let $S=\{a_1a_2\dots a_n:c(\{a_1,a_2,\dots,a_n\})=0\}$; then both $F_n^p\cup\{S\}$ and $F_n^p\cup\{P^{(n)}\setminus S\}$ have the f.i.p.\ so there are at least two ultrafilters containing $F_n^p$.\kraj

It is well-known that, under CH, there are Ramsey ultrafilters in $N^*$ (\cite{CN}, Theorem 9.19). Also, not all ultrafilters in $N^*$ are Ramsey. If $f:N\str P$ is any bijection, $\widetilde{f}$ maps Ramsey to Ramsey, and non-Ramsey to non-Ramsey ultrafilters, so (under CH) in $P^*$ there are also ultrafilters of both types. On the other hand, each Ramsey ultrafilter is a P-point and Shelah proved that it is consistent with ZFC that there are no P-points in $N^*$ (the proof can be found in \cite{Sh}).

Blass proved in \cite{B2} that (under CH) there is a non-Ramsey ultrafilter $p$ such that for every 3-coloring $c:[N]^2\str\{0,1,2\}$ there is $A\in p$ such that $|c[[A]^2]|\leq 2$. By a simple modification of the proof of Theorem \ref{ramsey} we get that, for such $p\in P^*$, there are exactly two ultrafilters $x\supseteq F_2^p$.

\begin{lm}\label{proizvod}
Let $x\in\overline{P^{(2)}}$.

(a) If there are disjoint $A,B\ps P$ such that $AB\in x$ then $x$ is divisible by at least two ultrafilters from $\overline{P}$.

(b) For any two distinct $p,q\in\overline{P}$, $x$ is divisible by $p$ and $q$ if and only if for every two disjoint $A\in p$, $B\in q$ holds $AB\in x$.
\end{lm}

\dokaz (a) If $AB\in x$ for disjoint $A,B\ps P$ then we can define functions $f_A$ and $f_B$ on $AB$ such that $f_A(ab)=a$ and $f_B(ab)=b$ (for $a\in A$, $b\in B$) so, by Lemma \ref{funkcija}(a), $\widetilde{f_A}(x)\widemid x$ and $\widetilde{f_B}(x)\widemid x$. Finally, $\widetilde{f_A}(x)\neq\widetilde{f_B}(x)$ because $A=f_A[AB]\in\widetilde{f_A}(x)$ and $B=f_B[AB]\in\widetilde{f_B}(x)$.\\

(b) First let $p\widemid x$ and $q\widemid x$. Let $A\in p\rest P$ and $B\in q\rest P$ be disjoint. Then $A\gstr\in p\cap\cU\subseteq x$ and $B\gstr\in q\cap\cU\subseteq x$, so $AB=A\gstr\cap B\gstr\cap P^{(2)}\in x$.

Now let $P=A\cup B$ be a partition of $P$ with $A\in p$, $B\in q$. As in (a) we define $f_A$ and $f_B$. Since $AB\in x$, $\widetilde{f_A}(x)$ and $\widetilde{f_B}(x)$ are ultrafilters below $x$. Assume $\widetilde{f_A}(x)\neq p$. Let $A_1\in\widetilde{f_A}(x)\setminus p$ be disjoint from $B$; by the first implication $A_1B\in x$. But $A_1B$ and $(P\setminus(A_1\cup B))B$ (belonging to $x$ by assumption) are disjoint so they can not both be in $x$, a contradiction. Thus $\widetilde{f_A}(x)=p$, and in the same way we prove $\widetilde{f_B}(x)=q$.\kraj

\begin{lm}\label{twodiv}
Let $p\in\overline{P}$. The ultrafilters such that their only proper divisors are 1 and $p$ are exactly $p^2$ and (for $p\in P^*$) ultrafilters containing $F_2^p$.
\end{lm}

\dokaz In Lemma \ref{kvadrati} we proved that the only proper divisors of $p^2$ are 1 and $p$. The proof for $x\supseteq F_2^p$ is similar. Now assume that $x$ is any ultrafilter such that 1 and $p$ are its only proper divisors. $x$ belongs either to $\overline{P}$, $\overline{P^2}$, $\overline{P^{(2)}}$ or to $\overline{(P\cup P^2\cup P^{(2)}\cup\{1\})^c}$.

In the first case $x$ has only one proper divisor, 1.

In the second case, if we let $p=\{A\ps P:A^2\in x\}$, it is easy to prove $x=p^2$.

Let $x\in\overline{P^{(2)}}$ and $G=\{A\subseteq P:A^{(2)}\in x\}$. $G$ is closed for finite intersections and sets in $G$ are nonempty, so $G$ has the f.i.p. If $A\cup B$ is any partition of $P$, by Lemma \ref{proizvod} $AB\notin x$, so it follows that $A^{(2)}\in x$ or $B^{(2)}\in x$ (because $P^{(2)}=A^{(2)}\cup B^{(2)}\cup AB$), hence $A\in G$ or $B\in G$. This means that $G$ is an ultrafilter on $P$ and $x\supseteq F_2^G$. Since $p$ is the only element of $\overline{P}$ below $x$, we have $x\supseteq F_2^p$.

Finally, if $x\notin\overline{P\cup P^2\cup P^{(2)}\cup\{1\}}$, we can define the function $g:N\setminus(P\cup\{1\})\str N$ by $g(a_1a_2\dots a_n)=a_1a_2$ (where $a_1\leq a_2\leq\dots\leq a_n$ are prime). Then $\widetilde{g}(x)\neq x$ is an element of $P^2\cup P^{(2)}$ below $x$, so $x$ has more than two proper divisors.\kraj

The following definition and lemmas will be used in the proof of Theorem \ref{2cFpn}.

\begin{de}
Let $X\subseteq N$ and let $d=\{X_k:k\in N\}$ be a partition of $X^{(n)}$. A set $A\subseteq X$ is $d$-thick if for all $m\in N$ and all finite partitions $A=A_1\cup A_2\cup\dots\cup A_m$ there is $i\leq m$ such that for every $k\in N$ $A_i^{(n)}\cap X_k\neq\emptyset$.
\end{de}

The condition of $d$-thickness strengthens the condition that $A^{(n)}$ intersects every $X_k$. The idea of this strengthening is to satisfy (b) of the lemma below.

\begin{lm}\label{thick}
Let $A\subseteq X$, $B\subseteq X$ and let $d=\{X_k:k\in N\}$ be a partition of $X^{(n)}$.

(a) If $A$ is $d$-thick and $A\subseteq B$, then $B$ is $d$-thick.

(b) If neither $A$ nor $B$ are $d$-thick, then $A\cup B$ is not $d$-thick.
\end{lm}

\dokaz (a) is obvious.\\

(b) By (a) we may assume without loss of generality that $A\cap B=\emptyset$. Let partitions $A=A_1\cup A_2\cup\dots\cup A_{m_A}$ and $B=B_1\cup B_2\cup\dots\cup B_{m_B}$ and $k_i,l_i\in N$ be such that $A_i^{(n)}\cap X_{k_i}=\emptyset$ for all $i\leq m_A$ and $B_i^{(n)}\cap X_{l_i}=\emptyset$ for all $i\leq m_B$. Then the partition $A\cup B=A_1\cup A_2\cup\dots\cup A_{m_A}\cup B_1\cup B_2\cup\dots\cup B_{m_B}$ and the same $k_i,l_i$ witness that $A\cup B$ is not $d$-thick.\kraj

\begin{lm}\label{thick2}
If, for each $n\geq 2$, $d_n=\{X_{n,k}:k\in N\}$ is a partition of $X^{(n)}$ such that for all $k\in N$
\begin{equation}\label{eqparticije2}
X_{n+1,k}\subseteq\{xa:x\in X_{n,k},a\in P\}
\end{equation}
and $A\subseteq X$ is not $d_{n_0}$-thick then, for all $n\geq n_0$, $A$ is not $d_n$-thick.
\end{lm}

\dokaz It suffices to prove the theorem for $n=n_0+1$. Let the partition $A=A_1\cup A_2\cup\dots\cup A_m$ and $k_i\in N$ be such that $A_i^{(n_0)}\cap X_{n_0,k_i}=\emptyset$ for all $i\leq m$. Then the same partition witnesses that $A$ is not $d_{n_0+1}$-thick: if $xa\in A_i^{(n_0+1)}\cap X_{n_0+1,k_i}$ (for $x\in X_{n_0,k}$, $a\in P$), then $x\in A_i^{(n_0)}\cap X_{n_0,k_i}$, a contradiction.\kraj

\begin{lm}\label{Pthick}
(a) There is a coloring $c:[N]^2\str N$ such that for every $k\in N$:
\parbox{11cm}{\begin{eqnarray*}
\mbox{for every arithmetic progression }s\mbox{ of length }2^{k}+1\nonumber\\
\mbox{ there are }a,b\in s\mbox{ such that }c(\{a,b\})=k.
\end{eqnarray*}}\hfill
\parbox{1cm}{\begin{eqnarray}\label{progr}\end{eqnarray}}

(b) There are partitions $d_n=\{X_{n,k}:k\in N\}$ of $P^{(n)}$ (for $n\geq 2$) such that $P$ is $d_n$-thick for every $n\geq 2$ and (\ref{eqparticije2}) holds.
\end{lm}

\dokaz (a) We define sets $A^0_i=\{i\}$ (for $i\in N$) and, by recursion on $n$, $A^n_i=A^{n-1}_{2i-1}\cup A^{n-1}_{2i}$. Note that $a\in A_{2i-1}^{n-1}$ are exactly numbers with residue $1\leq r\leq 2^{n-1}$ modulo $2^{n}$.

First we color the pairs of numbers in the same $A^1_i$: $c(\{2i-1,2i\})=1$. By recursion on $n$, pairs of numbers $a,b\in A^n_i$ such that $a\in A^{n-1}_{2i-1}$ and $b\in A^{n-1}_{2i}$ are colored according to the difference $b-a$: if $2^j\leq b-a<2^{j+1}$ then $c(\{a,b\})=n-j$.

Let us prove (\ref{progr}). Let $s=\{a_0+md:0\leq m\leq 2^{k}\}$ and $j$ is such that $2^j\leq d<2^{j+1}$. At least one $a_0+md\in s$ (for $0\leq m<2^{k}$) has residue $2^{k+j-1}-d<r\leq 2^{k+j-1}$ modulo $2^{k+j}$, which means that $a_0+md\in A^{k+j-1}_{2i-1}$ and $a_0+(m+1)d\in A^{k+j-1}_{2i}$ for some $i\in N$. Then $c(\{a_0+md,a_0+(m+1)d\})=(k+j)-j=k$.\\

(b) Let $c$ be the coloring of $[N]^2$ defined in (a). We define $c_n:[N]^n\str N$ (for $n\geq 2$) by $c_n(\{a_1,a_2,\dots,a_n\})=c(\{a_1,a_2\})$ for $a_1<a_2<\dots<a_n$. Then (\ref{progr}) implies that for every $n\geq 2$ and every $k\in N$

\parbox{10cm}{\begin{eqnarray*}
\mbox{for every arithmetic progression }s\mbox{ of length }2^{k}+n-1\nonumber\\
\mbox{ there are }a_1,a_2,\dots,a_n\in s\mbox{ such that }c_n(\{a_1,a_2,\dots,a_n\})=k.\label{progrn}
\end{eqnarray*}}\hfill
\parbox{0.5cm}{\begin{eqnarray}\label{progrn}\end{eqnarray}}
Now enumerate $P=\{p_n:n\in N\}$ in the increasing order. Define the partitions $d_n=\{X_{n,k}:k\in N\}$ of $P^{(n)}$: $X_{n,k}=\{p_{a_1}p_{a_2}\dots p_{a_n}:c_n(\{a_1,a_2,\dots,a_n\})=k\}$. Obviously, (\ref{eqparticije2}) holds.

To prove that $P$ is $d_n$-thick, let $P=A_1\cup A_2\cup\dots\cup A_m$ be any finite partition of $P$. If $B_i=\{a\in N:p_a\in A_i\}$, then $N=B_1\cup B_2\cup\dots\cup B_m$ is a partition of $N$. By the infinite Van der Waerden's theorem there is $i\leq m$ such that $B_i$ contains arithmetic progressions of any length. So for every $k\in N$, there is a progression $s$ in $B_i$ of length $2^{k}+n-1$. By (\ref{progrn}) there are $a_1,a_2,\dots,a_n\in s\subseteq B_i$ such that $c_n(\{a_1,a_2,\dots,a_n\})=k$ so $p_{a_1}p_{a_2}\dots p_{a_n}\in X_{n,k}\cap A_i^{(n)}$.\kraj

\begin{te}\label{2cFpn}
(CH) There is $p\in P^*$ such that for every $n\geq 2$ there are $2^{\goth c}$ ultrafilters $x\supseteq F_n^p$.
\end{te}

\dokaz Let, for $n\geq 2$, $d_n$ be the partition of $P^{(n)}$ given by Lemma \ref{Pthick}(b). By the Continuum Hypothesis we can enumerate all subsets of $P$ as $\langle S_\xi:\xi<\omega_1\rangle$; recall that $S_\xi'=P\setminus S_\xi$. By recursion on $\xi<\omega_1$ we define sets $A_\xi$ and families ${\cal F}_\xi$ such that:

(1$_\xi$) ${\cal F}_\xi$ is a countable family of infinite subsets of $P$, closed for finite intersections;

(2$_\xi$) ${\cal F}_\zeta\subseteq {\cal F}_\xi$ for $\zeta<\xi$;

(3$_\xi$) ${\cal F}_\xi=\bigcup_{\zeta<\xi}{\cal F}_\zeta$ for $\xi$ a limit ordinal;

(4$_\xi$) $A$ is $d_n$-thick for all $A\in {\cal F}_\xi$ and all $n\geq 2$;

(5$_\xi$) either $A_\xi=S_\xi$ or $A_\xi=S_\xi'$, and $A_\xi\in {\cal F}_{\xi+1}$.

First we let ${\cal F}_0=\{P\}$; by Lemma \ref{Pthick} $P$ is $d_n$-thick for all $n\geq 2$. Now for every $\xi<\omega_1$, assuming we have already defined ${\cal F}_\xi$ satisfying (1$_\xi$)-(4$_\xi$), we define $A_\xi$ and ${\cal F}_{\xi+1}$.

We first prove that for at least one of the possibilities $A_\xi=S_\xi$ and $A_\xi=S_\xi'$ all sets $A_\xi\cap A$ for $A\in F_\xi$ are $d_n$-thick for all $n\geq 2$. Assume not: then there are $A,B\in {\cal F}_\xi$ and $m,n\in N$ such that $S_\xi\cap A$ is not $d_m$-thick and $S_\xi'\cap B$ is not $d_n$-thick. If, say, $m\leq n$, by Lemma \ref{thick2} both those sets are not $d_n$-thick. By Lemma \ref{thick}(a), $S_\xi\cap(A\cap B)$ and $S_\xi'\cap(A\cap B)$ are not $d_n$-thick and, finally, by (b) of the same lemma, their union $P\cap(A\cap B)=A\cap B$ is not $d_n$-thick, which is impossible by (4$_\xi$) and (1$_\xi$).

Hence we define $A_\xi=S_\xi$ if all $S_\xi\cap A$ are $d_n$-dense for all $n\geq 2$ and all $A\in {\cal F}_\xi$, and $A_\xi=S_\xi'$ otherwise. Let ${\cal F}_{\xi+1}={\cal F}_\xi\cup\{A_\xi\cap A:A\in {\cal F}_\xi\}$. If $\xi$ is a limit ordinal, we define ${\cal F}_\xi$ as in (3$_\xi$). Clearly, all the properties (1$_\xi$)-(5$_\xi$) are now satisfied.

In the end, by (1$_\xi$) and (5$_\xi$) $p:=\bigcup_{\xi<\omega_1}{\cal F}_\xi$ is an ultrafilter on $P$. Let $n\geq 2$. For every $k\in N$ the family $F_n^p\cup\{X_{n,k}\}$ has the f.i.p.\ (since every $A\in p\rest P$ is $d_n$-thick, every $A^{(n)}\in F_n^p$ intersects $X_{n,k}$), so there are ultrafilters $x_k\supseteq F_n^p\cup\{X_{n,k}\}$. But $X_{n,k}$ (for $k\in N$) are disjoint, so all $x_k$ are distinct ultrafilters. By Lemma \ref{broj}(b) there are $2^{\goth c}$ ultrafilters $x\supseteq F_n^p$.\kraj


\section{Ultrafilters with two prime divisors}

\begin{lm}\label{atmosttwo}
No ultrafilter $x\in\overline{P^{(2)}}$ is divisible by more than two prime ultrafilters.
\end{lm}

\dokaz Assume the opposite, that $x$ is divisible by $p_1,p_2,p_3\in\overline{P}$. Let $P=A_1\cup A_2\cup A_3$ be a partition such that $A_1\in p_1$, $A_2\in p_2$ and $A_3\in p_3$. We consider three cases. $1^\circ$ If $A_1P\notin x$, then clearly $A_1\gstr\notin x$ (because $A_1\gstr\cap P^{(2)}=A_1P$), so $p_1\nwidemid x$. $2^\circ$ If $A_1(A_1\cup A_3)\in x$, then $A_2P\notin x$ so, as in $1^\circ$, $p_2\nwidemid x$. $3^\circ$ Finally, if $A_1A_2\in x$ then $A_3P\notin x$ so $p_3\nwidemid x$. In each case we reach a contradiction.\kraj

\begin{de}
$F_{1,1}^{p,q}=\{AB:A\in p\rest P,B\in q\rest P,A\cap B=\emptyset\}$ for $p,q\in\overline{P}$.
\end{de}

By Lemma \ref{proizvod} ultrafilters $x\supseteq F_{1,1}^{p,q}$ are exactly those ultrafilters in $P^{(2)}$ divisible by $p$ and $q$.

\begin{lm}\label{broj2}
For any distinct $p,q\in\overline{P}$:

(a) $p\cdot q\supseteq F_{1,1}^{p,q}$ and $q\cdot p\supseteq F_{1,1}^{p,q}$.

(b) there are either finitely many or $2^{\goth c}$ ultrafilters $r\supseteq F_{1,1}^{p,q}$.
\end{lm}

\dokaz (a) Let $A\in p$ and $B\in q$ be disjoint. Since $AB/a=B$ for $a\in A$, $\{a\in N:AB/a\in q\}\supseteq A\in p$ and so $AB\in p\cdot q$. Analogously $AB\in q\cdot p$.\\

(b) Analogous to the proof of Lemma \ref{broj}(b).\kraj

If $n\in N$ then $nq=qn$ for all $q\in\beta N$ and it is not hard to see that this is the only ultrafilter containing $x\supseteq F_{1,1}^{n,q}$: by Lemma 5.1 of \cite{So1} every ultrafilter divisible by $n$ must contain the set $nN=\{na:a\in N\}$ so $nP\in x$; but $F_{1,1}^{n,q}=\{n(B\setminus\{n\}):B\in q\}$ generates an ultrafilter on $nP$.

On the other hand, by \cite{HS}, Corollary 6.51, for every $p\in P^*$ there is $q\in P^*$ such that $pq\neq qp$, and we have at least two ultrafilters containing $F_{1,1}^{p,q}$. We will improve this in Theorem \ref{dvaiznadpq}.

\begin{te}
Let $p,q\in P^*$. If there is unique $x\supseteq F_{1,1}^{p,q}$ then both $p$ and $q$ are P-points.
\end{te}

\dokaz First let $X_n\in q$ for $n\in N$ and, without loss of generality, assume $X_0=P$ and $X_m\subseteq X_n$ for $m>n$. Let $Y=\{mn\in P^{(2)}:m>n\mbox{ and }m\in X_n\}$. Then, for every $n\in P$, $Y/n\supseteq X_n\setminus\{1,2,\dots,n\}\in q$ so $\{n\in N:Y/n\in q\}=P\in p$ and $Y\in p\cdot q$.

Since $F_{1,1}^{p,q}$ generates the unique ultrafilter, there is $AB\in F_{1,1}^{p,q}$ such that $AB\subseteq Y$ ($A\in p$, $B\in q$). To prove that $q$ is a P-point it suffices to show that $B$ is a pseudointersection of the sets $X_n$. For every $n\in N$ there is $a\in A$ such that $n\leq a$. If $b\in B$ is such that $b>a$, then $ab\in Y$ implies $b\in X_a\subseteq X_n$. Hence $B\setminus X_n$ is finite.

Since $p\cdot q$ and $q\cdot p$ both contain $F_{1,1}^{p,q}$, we have $p\cdot q=q\cdot p$, so by interchanging the roles of $p$ and $q$ we prove in the same way that $p$ is a P-point.\kraj

\begin{te}\label{dvaiznadpq}
For every $p\in P^*$ there is an ultrafilter $q\in P^*$ such that there are $2^{\goth c}$ ultrafilters $r\supseteq F_{1,1}^{p,q}$.
\end{te}

\dokaz Let $p\in P^*$ be given. Let $f:P\str \fun NP$ be such that, if $f(i)=f_i$, then $\langle f_i:i\in P\rangle$ is the sequence of all eventually constant functions $f_i:N\str P$ (i.e.\ such that there are $n_0\in N$ and $a\in P$ so that $f_i(n)=a$ for $n\geq n_0$), ordered in such way that for all $n\in N$: $f_i(n)\leq i$ for $i\in\{2,3\}$ and $f_i(n)<i$ for $i>3$. The set $D=\{f_i:i\in P\}$ is dense in the space $\fun NP$ (with the usual Tychonoff topology). For $m,n\in N$ and $A\in p$ let $U_n(A)=\{x\in\fun NP:x(n)\in A\}$ and $V_{m,n}=\{x\in\fun NP:x(m)\neq x(n)\}$. Then the family $\{U_n(A):n\in N\land A\in p\}\cup\{V_{m,n}:m\neq n\}$ has the uniform f.i.p.\ so the family $F=\{f^{-1}[U_n(A)]:n\in N\land A\in p\}\cup\{f^{-1}[V_{m,n}]:m\neq n\}$ also has the uniform f.i.p. (because $f$ is one-to-one, all $U_n(A)$ and $V_{m,n}$ are open and $D$ is dense in $\fun NP$). Hence there is an ultrafilter $q\in P^*$ containing $F$.

Now, for $n\in N$, let $g_n:P\str N$ be defined by $g_n(i)=if_i(n)$, and let $r_n=\widetilde{g_n}(q)$. Since $f_i(n)\neq i$ for all $i>3$ and all $n$, we conclude that $r_n\in\overline{P^{(2)}\cup\{4,9\}}$. Moreover, by Lemma \ref{funkcija}(b) all $r_n$ are divisible by $q$, hence they are nonprincipal and $r_n\in(P^{(2)})^*$.

Let $h_1,h_2:P^{(2)}\cup\{4,9\}\str P$ be defined by $h_1(ab)=a$ and $h_2(ab)=b$ for $a,b\in P$ and $a\leq b$. We prove that $\widetilde{h_1}(r_n)=p$ and $\widetilde{h_2}(r_n)=q$ for every $n\in N$. First, $\widetilde{h_2}(r_n)=\widetilde{h_2}\circ\widetilde{g_n}(q)=\widetilde{h_2\circ g_n}(q)=q$ (because $h_2\circ g_n$ is the identity function).

On the other hand, $h_1\circ g_n(i)=f_i(n)$ for all $i$. For any $A\in p$ and all $i\in f^{-1}[U_n(A)]$ we have $f_i(n)\in A$. Hence $h_1\circ g_n[f^{-1}[U_n(A)]]\subseteq A$ and so $A\in\widetilde{h_1\circ g_n}(q)$. Thus $\widetilde{h_1}(r_n)=\widetilde{h_1\circ g_n}(q)=p$.

By Lemma \ref{funkcija}(a) $p\widemid r_n$ and $q\widemid r_n$ for each $n\in N$, so $r_n\in F_{1,1}^{p,q}$. It remains to prove that all $r_n$ are distinct so, by Lemma \ref{broj2}(b), the set of ultrafilters $r\supseteq F_{1,1}^{p,q}$ will be of cardinality $2^{\goth c}$. Let $m<n$. We prove that the sets $g_m[f^{-1}[V_{m,n}]]$ and $g_n[f^{-1}[V_{m,n}]]$ are disjoint: if we assume that $g_m(i)=g_n(j)$ for some $i,j\in f^{-1}[V_{m,n}]$ then $if_i(m)=jf_j(n)$, so since $f_i(m)\leq i$, $f_j(n)\leq j$ and all of the numbers $i,j,f_i(m),f_j(n)$ are prime, we have $i=j$ and $f_i(m)=f_i(n)$, a contradiction with the fact $f_i\in V_{m,n}$. But $g_m[f^{-1}[V_{m,n}]]\in r_m$ and $g_n[f^{-1}[V_{m,n}]]\in r_n$, so $r_m\neq r_n$.\kraj

\section{The higher levels}

\begin{de}
We call ultrafilters of the form $p^k$ for some $p\in\overline{P}$ and $k\in N$ {\it basic}. Let $\cB$ be the set of all basic ultrafilters, and let $\cA$ be the set of all functions $\alpha:\cB\str N\cup\{0\}$ with finite support (i.e.\ such that $\{b\in\cB:\alpha(b)\neq 0\}$ is finite).
\end{de}

We will abuse notation and write $\alpha=\{(b_1,n_1),(b_2,n_2),\dots,(b_m,n_m)\}$ if $\alpha(b)=0$ for $b\notin\{b_1,b_2,\dots,b_m\}$ (allowing also some of the $n_i$ to be zeros).

\begin{de}
Let $\alpha=\{(p_1^{k_1},n_1),(p_2^{k_2},n_2),\dots,(p_m^{k_m},n_m)\}\in\cA$ ($p_i\in\overline{P}$). With $F_\alpha$ we denote the family of all sets
\begin{eqnarray*}\label{eqhigher}
(A_1^{k_1})^{(n_1)}(A_2^{k_2})^{(n_2)}\dots(A_m^{k_m})^{(n_m)}=\{\prod_{i=1}^m\prod_{j=1}^{n_i}a_{i,j}^{k_i} &:& a_{i,j}\in A_i\mbox{ for all }i,j\\
&&\land\mbox{ all }a_{i,j}\mbox{ are distinct}\}
\end{eqnarray*}
such that: (i) $A_i\in p_i\rest P$, (ii) $A_i=A_j$ if $p_i=p_j$ and $A_i\cap A_j=\emptyset$ otherwise.
\end{de}

\begin{ex}
$F_{1,1}^{p,q}=F_\alpha$ for $\alpha=\{(p,1),(q,1)\}$ and $F_2^p=F_\beta$ for $\beta=\{(p,2)\}$. If $\gamma=\{(p^2,1)\}$ then the set $F_\gamma=\{A^2:A\in p\rest P\}$ generates $p^2$.
\end{ex}

\begin{de}
If $\alpha=\{(p_1^{k_1},n_1),(p_2^{k_2},n_2),\dots,(p_m^{k_m},n_m)\}\in\cA$ ($p_i\in\overline{P}$), we denote $\sigma(\alpha)=\sum_{i=1}^m(k_in_i)$.
\end{de}

\begin{te}\label{nivoi}
The $n$-th level $\overline{L_n}$ (for $n\in N$) consists precisely of ultrafilters containing $F_\alpha$ for some $\alpha\in\cA$ such that $\sigma(\alpha)=n$.
\end{te}

\dokaz To avoid cumbersome notation we prove the theorem only for $n=4$. This special case contains essentially all the ideas needed for the proof in general.

Ultrafilters in the 4th level contain $L_4=\{a_1a_2a_3a_4:a_1,a_2,a_3,a_4\in P\}$. We partition $L_4$ as
$$L_4=P^4\cup P^3P\cup (P^2)^{(2)}\cup P^2P^{(2)}\cup P^{(4)}.$$
So for every ultrafilter $x\in\overline{L_4}$ we have 5 cases.

$1^\circ$ $x\ni P^4=\{a^4:a\in P\}$. Then $p:=\{A\ps P:A^4\in x\}$ is an ultrafilter in $\overline{P}$ and it is easy to see that $x=p^4$, so $x\supseteq F_\alpha$ for $\alpha=\{(p^4,1)\}$.

$2^\circ$ $x\ni P^3P=\{a^3b:a,b\in P,a\neq b\}$. We define functions $f_1(a^3b)=a^3$ and $f_2(a^3b)=b$ for $a^3b\in P^3P$. Let $p^3={\widetilde f_1}(x)$ (clearly, ${\widetilde f_1}(x)\in\overline{P^3}$) and $q={\widetilde f_2}(x)$. We prove that $x\supseteq F_\alpha$ for $\alpha=\{(p^3,1),(q,1)\}$. Let $A^3B\in F_\alpha$ (it may be that $p=q$ - then $A=B$, or $p\neq q$ and $A\cap B=\emptyset$). Then $A^3B=f_1^{-1}[A^3]\cap f_2^{-1}[B]\in x$.

$3^\circ$ $x\ni (P^2)^{(2)}=\{a^2b^2:a,b\in P,a\neq b\}$. We define functions $f_1(a^2b^2)=a^2$ and $f_2(a^2b^2)=b^2$ for $a^2b^2\in (P^2)^{(2)}$ and $a<b$. Let $p^2={\widetilde f_1}(x)$ and $q^2={\widetilde f_2}(x)$. If $p\neq q$ then $x\supseteq F_\alpha$ for $\alpha=\{(p^2,1),(q^2,1)\}$: for each disjoint $A\in p\rest P$ and $B\in q\rest P$, $A^2B^2\supseteq f_1^{-1}[A^2]\cap f_2^{-1}[B^2]\in x$. Otherwise, $x\supseteq F_\alpha$ for $\alpha=\{(p^2,2)\}$: for each $A\in p$, $(A^2)^{(2)}=f_1^{-1}[A^2]\cap f_2^{-1}[A^2]\in x$.

$4^\circ$ $x\ni P^2P^{(2)}=\{a^2bc:a,b,c\in P,a\neq b\neq c\neq a\}$. We define functions $f_1(a^2bc)=a^2$, $f_2(a^2bc)=b$ and $f_3(a^2bc)=c$ for $a^2bc\in P^2P^{(2)}$ and $b<c$. Let $p^2={\widetilde f_1}(x)$, $q={\widetilde f_2}(x)$ and $r={\widetilde f_3}(x)$. If $q=r$ then $x\supseteq F_\alpha$ for $\alpha=\{(p^2,1),(q,2)\}$, otherwise $x\supseteq F_\alpha$ for $\alpha=\{(p^2,1),(q,1),(r,1)\}$.

$5^\circ$ $x\ni P^{(4)}=\{abcd:a,b,c,d\in P,\;a,b,c,d\mbox{ all distinct}\}$. Analogously to previous cases, $x\supseteq F_\alpha$ for one of the following: $\alpha=\{(p,4)\}$, $\alpha=\{(p,3),(q,1)\}$, $\alpha=\{(p,2),(q,1),(r,1)\}$, $\alpha=\{(p,2),(q,2)\}$ or $\alpha=\{(p,1),(q,1),(r,1),(s,1)\}$ for some $p,q,r,s\in\overline{P}$.\kraj

\begin{de}
For every $\alpha\in\cA$ and every $p\in\overline{P}$ let $\alpha\rest p=\langle\alpha(p^k):k\in N\rangle$. Clearly, all such sequences have finitely many non-zero elements. If $\vec{x}=\langle x_k:k\in N\rangle$ and $\vec{y}=\langle y_k:k\in N\rangle$ are two such sequences in $N\cup\{0\}$ we say that $\vec{y}$ {\it dominates} $\vec{x}$ if for every $m\in N$, $\sum_{k\geq m}x_k\leq\sum_{k\geq m}y_k$.

We define an order on $\cA$ as follows: $\alpha\leq\beta$ if for every $p\in\overline{P}$ $\beta\rest p$ dominates $\alpha\rest p$.
\end{de}

\begin{ex}
If $\alpha=\{(p,2)\}$ and $\beta=\{(p,1),(p^2,1)\}$ (for some $p\in\overline{P}$), then $\alpha\rest p=\langle 2,0,0,\dots\rangle$ and $\beta\rest p=\langle 1,1,0,0,\dots\rangle$ so $\beta\rest p$ dominates $\alpha\rest p$ and $\alpha\leq\beta$. 

But if $\alpha=\{(p,2),(p^2,2),(q,2)\}$ and $\beta=\{(p,1),(p^3,3),(q^2,1)\}$, then $\beta\rest p=\langle 1,0,3,0,\dots\rangle$ dominates $\alpha\rest p=\langle 2,2,0,0,\dots\rangle$, but $\beta\rest q=\langle 0,1,0,0,\dots\rangle$ does not dominate $\alpha\rest q=\langle 2,0,0,\dots\rangle$ so $\alpha\leq\beta$ does not hold.
\end{ex}

It is not hard to see that $\alpha\leq\beta$ implies $\sigma(\alpha)\leq\sigma(\beta)$.

\begin{te}\label{poredak}
Let $\alpha,\beta\in\cA$.

(a) If $x\supseteq F_\alpha$, $y\supseteq F_\beta$ and $x\widemid y$, then $\alpha\leq\beta$.

(b) There is no $x\supseteq F_\alpha\cup F_\beta$ for $\alpha\neq\beta$.
\end{te}

\dokaz (a) Assume that $\alpha\leq\beta$ does not hold. This means that, for some prime $p$, $\beta\rest p=\langle y_k:k\in N\rangle$ does not dominate $\alpha\rest p=\langle x_k:k\in N\rangle$, i.e.\ there is $m\in N$ such that
$\sum_{k\geq m}x_k>\sum_{k\geq m}y_k.$ Let $n=\max\{k:x_k>0\lor y_k>0\}$; then
\begin{equation}\label{eqsigma}
u:=\sum_{k=m}^nx_k>v:=\sum_{k=m}^ny_k.
\end{equation}
For any $A\in p$, if we denote $B=((A^m)^{(x_m)}(A^{m+1})^{(x_{m+1})}\dots(A^n)^{(x_n)})\gstr$, we have $B\in x$, because $B\cap L_{\sigma(\alpha)}$ contains a set in $F_\alpha$. On the other hand, $B\notin y$, because $B$ is disjoint from any set in $F_\beta$: every element of $B$ has a divisor of the form $a_{m,1}^m\dots a_{m,x_m}^ma_{m+1,1}^{m+1}\dots a_{n,x_n}^n$ with $u$ prime factors from $A$ to powers $\geq m$ and, because of (\ref{eqsigma}), no element from a set in $F_\beta$ does (they all have exactly $v$ prime factors from $A$ to powers $\geq m$). Thus, $B\in x\cap\cU\setminus y$. This is a contradiction with $x\widemid y$.\\

(b) Since $x\widemid x$, by (a) $x\supseteq F_\alpha\cup F_\beta$ would imply $\alpha\leq\beta$ and $\beta\leq\alpha$. However, the relation $\leq$ on $\cA$ is clearly antisymmetric, so $\alpha=\beta$.\kraj

\begin{ex}\label{example}
The reverse of Theorem \ref{poredak}(a) does not hold: we will construct ultrafilters $x\supseteq F_\alpha$ and $y\supseteq F_\beta$ such that $\alpha\leq\beta$ but $x\nwidemid y$.

If $p\in P^*$ is not Ramsey, by Theorem \ref{ramsey} there is $A\subseteq P^{(2)}$ such that neither $A$ nor $P^{(2)}\setminus A$ contain a set in $F_2^p$. Exactly one of these two sets is in $p\cdot p$, say $A\in p\cdot p$, i.e.\ $\{a\in P:A/a\in p\}\in p$. Then for every $X\in p$ we can choose $a\in X$ such that $A/a\in p$. If $bc\notin A$ for all distinct $b,c\in X\cap A/a$, then $(X\cap A/a)^{(2)}\subseteq P^{(2)}\setminus A$, a contradiction with our choice of $A$. So for every $X\in p$ there are $a,b,c\in X$ such that $ab,ac,bc\in A$. This means that, if we denote $S=\{abc:ab,ac,bc\in A\}$, the family $F_3^p\cup\{S\}$ has the f.i.p. Hence we can find ultrafilters $x\supseteq F_2^p\cup\{P^{(2)}\setminus A\}$ and $y\supseteq F_3^p\cup\{S\}$, so $S\dstr\in(y\cap\cV)\setminus x$ and thus $x\nwidemid y$.
\end{ex}

\begin{te}\label{ispodiznad}
Let $\alpha,\beta\in\cA$ be such that $\alpha\leq\beta$.

(a) If $x\supseteq F_\alpha$ then there is at least one $y\supseteq F_\beta$ such that $x\widemid y$.

(b) If $y\supseteq F_\beta$ then there is at least one $x\supseteq F_\alpha$ such that $x\widemid y$.
\end{te}

\dokaz (a) Let $\alpha=\{(b_1,n_1),(b_2,n_2),\dots,(b_m,n_m)\}$ and $\beta=\{(b_1,n_1'),(b_2,n_2')$, $\dots,(b_m,n_m')\}$. We prove that the family $F_\beta\cup(x\cap\cU)$ has the f.i.p. Since every intersection of finitely many elements of $F_\beta$ contains an element of $F_\beta$ and $x\cap\cU$ is closed for finite intersections, it suffices to prove that every set $B'=A_1^{(n_1')}A_2^{(n_2')}\dots A_m^{(n_m')}\in F_\beta$ intersects every $X\in x\cap\cU$. But $X$ intersects $B=A_1^{(n_1)}A_2^{(n_2)}\dots A_m^{(n_m)}\in F_\alpha$, say $l\in B\cap X$. It suffices to show that $B'$ contains numbers divisible by $l$.

Let $p$ be any prime ultrafilter which is a divisor of $x$. Let $l=a_1^{k_1}a_2^{k_2}\dots a_r^{k_r}b$ be the factorization of $l$ such that $a_1^{k_1},a_2^{k_2},\dots$ (where $a_i\neq a_j$ for $i\neq j$) are elements of those sets $A_1,A_2,\dots$ that belong to a power of $p$, ordered so that $k_1\geq k_2\geq\dots$, and $b$ is the product of the other factors of $l$. Let $l'=(a_1')^{k_1'}(a_2')^{k_2'}\dots (a_s')^{k_s'}b'$ be the factorization of any element $l'\in B'$ obtained in the same way. Then the fact that $\beta\rest p$ dominates $\alpha\rest p$ implies that $s\geq r$ and $k_i'\geq k_i$ for $i\leq r$. We let $l_p=a_1^{k_1'}a_2^{k_2'}\dots a_r^{k_r'}(a_{r+1}')^{k_{r+1}'}\dots(a_s')^{k_s'}$; this element is clearly divisible by $a_1^{k_1}a_2^{k_2}\dots a_r^{k_r}$. In the same way we construct, for every prime $q\widemid x$, a number $l_q$. The product of all such $l_q$ will be the desired element divisible by $l$.\\

(b) Analogously to (a), we can prove that $F_\alpha\cup(y\cap\cV)$ has the f.i.p.\kraj

Theorem \ref{ramsey} shows that the ultrafilter $y$ from Theorem \ref{ispodiznad}(a) need not be unique. The next example shows the same for $x$ from Theorem \ref{ispodiznad}(b).

\begin{ex}
As in Example \ref{example}, let $p\in P^*$ be non-Ramsey and let $A\subseteq P^{(2)}$ be such that both $F_2^p\cup\{A\}$ and $F_2^p\cup\{A"\}$ (where $A"=P^{(2)}\setminus A$) have the f.i.p. Then $F_4^p\cup\{AA"\}$ also has the f.i.p.: if $B^{(4)}\in F_4^p$ then $B^{(2)}$ intersects both $A$ and $A"$, say $a=b_1b_2\in B^{(2)}\cap A$ and $a"=b_1"b_2"\in B^{(2)}\cap A"$ (and without loss of generality $b_1",b_2"$ are different from $b_1,b_2$). But then $aa"\in B^{(4)}\cap AA"$. Hence there is $y\supseteq F_4^p\cup\{AA"\}$ and, if we define $f_1(aa")=a$ and $f_2(aa")=a"$ for $a\in A$, $a"\in A"$, then $\widetilde{f_1}(y)\supseteq F_2^p\cup\{A\}$ and $\widetilde{f_2}(y)\supseteq F_2^p\cup\{A"\}$ are both divisors of $y$.
\end{ex}

\begin{te}
If $x\supseteq F_\alpha$ and $y\supseteq F_\beta$, then $x\cdot y\supseteq F_{\alpha+\beta}$, where $\alpha+\beta=\{(b,n_b+n_b'):(b,n_b)\in x,(b,n_b')\in y\}$.
\end{te}

\dokaz First, the theorem clearly holds if one of the ultrafilters $x$, $y$ is in $N$. If $x=kx'$ and $y=ly'$ for some $k,l\in N$, $x',y'\in N^*$, then $x\cdot y=kl(x'\cdot y')$, so we can assume without loss of generality that $x$ and $y$ are not divisible by any elements of $N$.

Now let $\alpha=\{(b_1,n_1),(b_2,n_2),\dots,(b_m,n_m)\}$ and $\beta=\{(b_1,n_1'),(b_2,n_2'),\dots$, $(b_m,n_m')\}$. Then the sets in $F_{\alpha+\beta}$ are of the form $A_{b_1}^{(n_1+n_1')}A_{b_2}^{(n_2+n_2')}\dots A_{b_m}^{(n_m+n_m')}$ (where $A_p\in p$ for prime divisors $p$ of $x$ or $y$, $A_p\cap A_q=\emptyset$ for $p\neq q$ and, for basic divisors of the form $p^k$, $A_{p^k}=A_p^k$). But if $a=a_{1,1}a_{1,2}\dots a_{1,n_1}a_{2,1}\dots a_{m,n_m}\in A_{b_1}^{(n_1)}A_{b_2}^{(n_2)}\dots A_{b_m}^{(n_m)}$, then
$$(A_{b_1}^{(n_1+n_1')}A_{b_2}^{(n_2+n_2')}\dots A_{b_m}^{(n_m+n_m')})/a\supseteq B_{b_1}^{(n_1')}\dots B_{b_m}^{(n_m')},$$
where $B_{b_k}=A_{b_k}\setminus\{a_{1,1},a_{1,2},\dots a_{m,n_m}\}$. For every $a\in A_{b_1}^{(n_1)}A_{b_2}^{(n_2)}\dots A_{b_m}^{(n_m)}$ we have $B_{b_1}^{(n_1')}\dots B_{b_m}^{(n_m')}\in y$ so
$$\{a\in N:(A_{b_1}^{(n_1+n_1')}A_{b_2}^{(n_2+n_2')}\dots A_{b_m}^{(n_m+n_m')})/a\in y\}\supseteq A_{b_1}^{(n_1)}A_{b_2}^{(n_2)}\dots A_{b_m}^{(n_m)}\in x,$$
i.e.\ $A_{b_1}^{(n_1+n_1')}A_{b_2}^{(n_2+n_2')}\dots A_{b_m}^{(n_m+n_m')}\in x\cdot y$.\kraj

In particular, if $x\in\overline{L_m}$ and $y\in\overline{L_n}$ then $xy\in\overline{L_{m+n}}$. A corollary of this theorem is that levels $\overline{L_n}$ contain no idempotents (ultrafilters such that $x\cdot x=x$), since if $x\in\overline{L_n}$, then $x\cdot x\in\overline{L_{2n}}$. Another corollary is that, if $p\in\overline{P}$ and $x$ and $y$ are on finite levels, then $p\widemid x\cdot y$ implies $p\widemid x$ or $p\widemid y$, which (partly) justifies our calling such ultrafilters prime.

\begin{lm}
$|[x]_\sim|=1$ for every $n\in N$ and all ultrafilters $x\in\overline{L_n}$.
\end{lm}

\dokaz Assume the opposite, that $y=_\sim x$. If $x\in F_\alpha$ and $y\in F_\beta$, then Theorem \ref{poredak}(a) implies that $\alpha=\beta$, so $y\in\overline{L_n}$ as well. Let $A\ps L_n$ be such that $A\in y$ and $L_n\setminus A\in x$. Then $A\gstr\in y\cap\cU$ and $A\gstr\notin x$, a contradiction.\kraj

All of the proofs in the next lemma are analogous to Lemma \ref{broj}(b).

\begin{lm}
(a) For every $\alpha\in\cA$ there are either finitely many or $2^{\goth c}$ ultrafilters containing $F_\alpha$.

(b) For every $x\in\beta N$ there are either finitely many or $2^{\goth c}$ ultrafilters above $x$, and either finitely many or $2^{\goth c}$ ultrafilters below $x$.

(c) For every $x\supseteq F_\alpha$ and all $\beta\leq\alpha\leq\gamma$ there are either finitely many or $2^{\goth c}$ ultrafilters above $x$ in $F_\gamma$, and either finitely many or $2^{\goth c}$ ultrafilters below $x$ in $F_\beta$.
\end{lm}

\section{Above all finite levels}\label{above}

\begin{fa}\label{greatest}
(\cite{So2}, Theorem 4.1) There is the $\widemid$-greatest class (consisting of ultrafilters divisible by all other ultrafilters).
\end{fa}

\begin{lm}
Let $\langle x_n:n\in N\rangle$ be a sequence of ultrafilters such that $x_n\in\overline{L_n}$ and $x_m\widemid x_n$ for $m<n$. Then there is an ultrafilter divisible by all $x_n$ and not divisible by any ultrafilter which is not below any $x_n$.
\end{lm}

\dokaz Let $F_1=\{A\in\cU:A\in x_n\mbox{ for some }n\in N\}$ and $F_2=\{B\in\cV:B^c\notin x_n\mbox{ for all }n\in N\}$. We prove that $F_1\cup F_2$ has the f.i.p. First, both $F_1$ and $F_2$ are closed for finite intersections. So let $A\in F_1$ and $B\in F_2$. There is $n\in N$ such that $A\in x_n$. Since also $B^c\notin x_n$, we have $A\cap B\in x_n$, so $A\cap B\neq\emptyset$.\kraj

In particular, for every countable set $S\subseteq\overline{P}$ there is an ultrafilter divisible by all $p\in S$ but not divisible by any prime $p\in\overline{P}\setminus S$. Since there are $2^{2^{\goth c}}$ subsets of $\overline{P}$, this can not hold for all uncountable $S\subseteq\overline{P}$.

Let us also note that it was shown in \cite{So1} that $\widemid$ is not antisymmetric, so there are ultrafilters above all finite levels such that $|[x]_\sim|>1$.

\section{Final remarks}

The ultrafilters containing $F_2^p$ and $F_{1,1}^{p,q}$ bear similarities with ultrafilters that are preimages under the natural map from $\beta(N\times N)$ to $\beta N\times\beta N$. Such ultrafilters were investigated, among other papers, in \cite{B1}, \cite{H1} and \cite{BM}. Hence some of the proofs in this paper are modifications of ideas from these papers (most of which can also be found in Chapter 16 of \cite{CN}). Since some of these modifications were not quite trivial, and for the sake of completeness, we decided to include all the proofs here. Also, the proof of Lemma \ref{Pthick}(b) is more general, different and (hopefully) more intuitive then in \cite{H1}.



\footnotesize
\begin{thebibliography}{99}
\bibitem{B1}
       A.\ Blass,
       Orderings on ultrafilters.
       PhD thesis, 1970.
\bibitem{B2}
       A.\ Blass,
       Ultrafilter mappings and their Dedekind cuts,
       Trans.\ Amer.\ Math.\ Soc.\ 188 (1974), No.2, 327-340.
\bibitem{BM}
       A.\ Blass, G.\ Moche,
       Finite preimages under the natural map from $\beta(N\times N)$ to $\beta N\times\beta N$,
       Topology Proc.\ 26 (2002), 407-432.
\bibitem{CN}
       W.\ W.\ Comfort, S.\ Negrepontis,
       The theory of ultrafilters.
       Springer-Verlag, 1974.
\bibitem{H1}
       N.\ Hindman,
       Preimages of points undeer the natural map from $\beta(N\times N)$ to $\beta N\times\beta N$,
       Proc.\ Amer.\ Math.\ Soc.\ 37 (1973), 603-608.
\bibitem{HS}
       N.\ Hindman, D.\ Strauss:
       Algebra in the Stone-\v Cech compactification, theory and applications.
       2nd revised and extended edition, De Gruyter, 2012.
\bibitem{Sh}
       S.\ Shelah,
       Proper and improper forcing,
       Perspectives in Mathematical Logic, Springer, 1998.
\bibitem{So1}
       B.\ \v Sobot,
       Divisibility in the Stone-\v Cech compactification,
       Rep.\ Math.\ Logic 50 (2015), 53-66.
\bibitem{So2}
       B.\ \v Sobot,
       Divisibility orders in the Stone-\v Cech compactification,
       submitted.
\bibitem{W}
       R.\ Walker,
       The \v Cech-Stone compactification.
       Springer, 1974.
\end{thebibliography}
\end{document}